\documentclass[a4paper,12pt]{amsart}

\usepackage{amssymb}
\usepackage{mhequ}
\renewcommand{\div}{\mbox{div}}

\usepackage{cite}

\DeclareFontFamily{OT1}{rsfs}{}
\DeclareFontShape{OT1}{rsfs}{m}{n}{ <-7> rsfs5 <7-10> rsfs7 <10-> rsfs10}{}
\DeclareMathAlphabet{\mathscr}{OT1}{rsfs}{m}{n}

\newcommand{\bel}[1]{\begin{equation}\label{#1}}
\newcommand{\beal}[1]{\begin{eqnarray}\label{#1}}
\newcommand{\beadl}[1]{\begin{deqarr}\label{#1}}
\newcommand{\eeadl}[1]{\arrlabel{#1}\end{deqarr}}
\newcommand{\eeal}[1]{\label{#1}\end{eqnarray}}
\newcommand{\eead}[1]{\end{deqarr}}
\newcommand{\eea}{\end{eqnarray}}
\newcommand{\eeaa}{\end{eqnarray*}}

\newcommand{\be}{\begin{equation}}
\newcommand{\ee}{\end{equation}}

\DeclareFontFamily{OT1}{rsfs}{}
\DeclareFontShape{OT1}{rsfs}{m}{n}{ <-7> rsfs5 <7-10> rsfs7 <10->
rsfs10}{} \DeclareMathAlphabet{\mycal}{OT1}{rsfs}{m}{n}

\newcounter{mnotecount}[section]

\newcommand{\rmnote}[1]{}

\newcommand{\Ric}{\operatorname{Ric}}

\def\mysavedown#1{\edef\mysubs{\mysubs#1}}
\def\mysaveup#1{\edef\mysups{\mysups#1}}
\def\mydown#1{{\mytensor}_{\vphantom{\mysubs}#1}}
\def\myup#1{{\mytensor}^{\vphantom{\mysups}#1}}
\def\tensor#1#2{
  #1
  \def\mytensor{\vphantom{#1}}
  \def\mysubs{\relax}
  \def\mysups{\relax}
  \let\down=\mysavedown
  \let\up=\mysaveup
  #2
  \let\down=\mydown
  \let\up=\myup
  #2
  }

\newcommand{\Tr}{\operatorname{Tr}}

\newcommand{\R}{\mathbb R}

\renewcommand{\div}{\operatorname{div}}

\renewcommand{\epsilon}{\varepsilon}
\renewcommand{\hat}{\widehat}

\def\crn#1#2{{\vcenter{\vbox{
        \hbox{\kern#2pt \vrule width.#2pt height#1pt
           }
          \hrule height.#2pt}}}}

\renewcommand{\hbar}{{\overline h}}

\newcommand{\pre}[2]{{{\vphantom{#2}}^{#1}}\kern-.2ex{#2}}

\theoremstyle{plain}

\theoremstyle{definition}

\numberwithin{equation}{section}

\date{June 2, 2016}

\begin{document}
\title[Conformally covariants  Lichnerowicz-York systems ]
{ Conformally covariant  parameterizations  for relativistic initial data}
\author[E. Delay]{Erwann Delay}
\address{Erwann Delay,
Laboratoire de Mathématiques d'Avignon (EA 2151), F-84018 Avignon, France}
\email{Erwann.Delay@univ-avignon.fr}
\urladdr{http://www.univ-avignon.fr/fr/recherche/annuaire-chercheurs\newline$\mbox{ }$\hspace{3cm}
/membrestruc/personnel/delay-erwann-1.html}

\begin{abstract}
 We revisit the Lichnerowicz-York method, and an alternative method of York, in order to obtain some conformally covariant systems.  This  type of parameterization is certainly more natural for non constant mean curvature initial data.
\end{abstract}

\maketitle

\\
\newline
\section{Introduction}
On a smooth manifold $M$ of dimension $n$, given a riemannian metric $g$, we denote by
 $R_g$ its scalar curvature and  $\nabla$ its Levi-Civita connection. If $h$ is a  symmetric covariant two tensor field, we define its divergence, as the 1-form given by
$$
(\div_g h)_i=-\nabla^kh_{ki}.
$$ 
The vacuum initial data for the  relativistic Einstein equations are given by a riemannian metric $\hat g$
and  a symmetric two tensor field $\hat K$, satisfying the constraint equations 
$$
(C)\left\{\begin{array}{l}R_{\hat g}-|\hat K|^2_{\hat g}+(\Tr_{\hat g}\hat K)^2=0\\
\div_{\hat g}\hat K+d(\Tr_{\hat g}\hat K)=0\end{array}\right.\;.
$$ 
This system is highly under-determined because it contains 
$(n+1)$ equations for $n(n+1)$
unknowns . It is natural to fix some of the unknowns, and to look for the  $(n+1)$ remaining ones.\\

In the usual conformal parameterization,  appears the 1944  Lichnerowicz equation 
\cite{Lichnerowicz:44} together with   the  1973 York decomposition
\cite{York:decompo}. It was actively  studied by many authors, on compact or non compact manifolds (asymptotic to some models :
euclidean, hyperbolic, cylindrical,...).

For a glimpse at the subject, we just mention the recent  paper
\cite{MaxwellConfetCTS} and its references. Otherwise, too many articles should be quoted.

We only recall that  after the  	
approaches realized for constant mean curvature $\tau$, like  \cite{IsenbergCMC} for instance,
several attempts were made in order to allow for variable mean curvature. 

The classical method  (historically method A, see \cite{BartnikIsenberg2004} section 4.1 for instance,
also called conformal TT method \cite{CookRelatNum2002})  starts with a given metric $g$, together with a trace free and divergence free symmetric two tensor $\sigma$ (a TT-tensor), and a function  $\tau$. One looks for the solutions of $(C)$ of the form
$$
\hat g=\phi^{N-2}g\;,\;\; \hat K=\frac{\tau}n\hat g+\phi^{-2}(\sigma+\mathring {\mathcal L}_g W),
\;\;\;(P)$$
where $N=\frac{2n}{n-2}$, and the unknowns are a function  $\phi>0$ and a one form $W$ and where
$$
(\mathring{\mathcal L}_g W)_{ij} =\nabla_iW_j+\nabla_jW_i-\frac2n\nabla^kW_k \;g_{ij}.
$$ Setting
$$
w=|\sigma+\mathring {\mathcal L}_g W|_g,
$$
we infer from $(C)$ and $(P)$ the coupled system
$$
(S)\left\{\begin{array}{ll}
P_{g,w}\phi:=\frac{4(n-1)}{n-2}\nabla^*\nabla\phi+R_g\phi-w^2\phi^{-N-1}=\frac{1-n}{n}\tau^2\phi^{N-1} &(L)\\
\div_g\mathring {\mathcal L}_gW=\frac{1-n}{n}\phi^{N}d\tau &(V)
\end{array}
\right.
$$
The Lichnerowicz equation $(L)$ has a covariant conformal property
(see section \ref{cc} for a precise definition). 
Indeed, if  $\phi$ is a solution of $(L)$, and $\varphi$ is any positive function, we may define
$$
\tilde g=\varphi^{N-2}g\;,\;\;\tilde w=\varphi^{-N}w\;,\;\; \tilde \phi=\varphi^{-1}\phi.
$$
Doing so, we find
$$
P_{\tilde g,\tilde w}\tilde \phi=\varphi^{1-N}P_{g,w}\phi=\frac{1-n}{n}\tau^2\tilde\phi^{N-1}.
$$

In contrast to $(L)$, the vectorial equation $(V)$ does not possess such a property. Moreover,
the conformal  transformation of $w$ used here for  $(L)$ does not correspond to a natural transformation of $\sigma$ and $W$.

However we note that  the operators appearing in $(V)$ are separately conformally covariant.
Specifically
$$
\div_{\tilde g} h=\varphi^{-N}\div_g(\varphi^2 h)\;,\;\;\mathring {\mathcal L}_{\tilde g}X=\varphi^{N-2}\mathring {\mathcal L}_g(\varphi^{2-N}X).\;\;\;(\mathcal C)
$$
But  the vector Laplacian  $\div\circ\mathring{\mathcal L}$ is not so.\\

Finally, we recall the   York decomposition \cite{York:decompo} valid for instance if $M$ is compact and  $g$ has no
conformal Killing fields  (i.e. ker$\;\mathring{\mathcal L}_g$ is trivial). Any  covariant symmetric trace free  two tensor field $h$ splits in a unique way as 
$$
 h=\sigma+\mathring{\mathcal L}_gW,\;\;\;(Y)
$$
where $\sigma$ is a TT-tensor and   $W$ a   1-form.\\

We define  in the section 
  \ref{paramCC} a  parameterization (method B of York)  derived from the  York decomposition relative to
  $\hat g$ instead of $g$.
This parameterization  gives rise to a  conformally covariant system.

The latter one leads to  a new vectorial laplacian, which is self adjoint for 
a weighted measure. This laplacian  will be studied in section \ref{equaV}.

In section \ref{methodBbis} we further discuss some aspects of the method B.
Amazingly, this method seems widely ignored in the literature neither pure nor numerical (see however \cite{CookRelatNum2002}).

In section \ref{autreparam} we  propose further parameterizations  (including conformally covariant ones) and compare them.

Finally, in section \ref{conj},  we give some comments related to the continuation of this work.\\

\noindent{\sc Acknowledgements}: I am grateful to Ph. Delano\"e for his multiples 
comments and to J. Isenberg for pointing out to me that J. York \cite{York:decompo} had already discovered the method B, as independently proposed in the earlier version of this paper.
I also thank P. Chru\'sciel  and  R.  Gicquaud for their comments.

\section{A first conformally covariant parameterization}\label{paramCC}
We recall here the method B of York (see section 4.1 of \cite{BartnikIsenberg2004} or the original paper \cite{York:decompo} page 461), also called physical TT method (see 
\cite{CookRelatNum2002} for instance).
Getting back to  $(C)$ and now using the York decomposition relative to  $\hat g$, namely
$$
\hat K-\frac{\tau}n\hat g=\hat\sigma+\mathring {\mathcal L}_{\hat g}\hat W,
$$
we find  that $(C)$ is equivalent to
$$
(C')\left\{\begin{array}{l}
\displaystyle{R_{\hat g}-|\hat\sigma+\mathring {\mathcal L}_{\hat g}\hat W|^2_{\hat g}=\frac{1-n}n\tau^2}\\
\displaystyle{\div_{\hat g}\mathring {\mathcal L}_{\hat g}\hat W=\frac{1-n}nd\tau}\\
\displaystyle{\div_{\hat g}\hat\sigma=0}
\end{array}\right.
$$
If we are looking for a solution in a conformal class $\hat g=\phi^{N-2}g$,
from $(\mathcal C)$, it is natural to introduce $\hat\sigma=\phi^{-2}\sigma $,
where $\sigma$ is a TT-tensor.  The third equation is then automatically satisfied.
Still using   $(\mathcal C)$, and expressing its  second equation in terms  of $g=
\phi^{2-N} \hat g$,
we are prompted to set   $\hat W:=\phi^{N-2} W$.

Sticking to the same fixed  $g,\sigma,\tau$, we can thus  parametrize the solutions of the constraint $(C)$ by
$$
(P')\;\left\{\begin{array}{l}
\hat g=\phi^{N-2}g\;,\\
 \hat K=\frac{\tau}n\hat g+\phi^{-2}(\sigma+\phi^N\mathring {\mathcal L}_g W)
\end{array}\right.
$$
With this parameterization in  $(C)$,  setting
$$
\omega=\omega(\sigma,\phi,W,g):=|\sigma+\phi^N\mathring {\mathcal L}_g W|_g,
$$
we obtain the new system
$$
(S')\left\{\begin{array}{ll}
P_{g,\omega}\phi=\frac{1-n}{n}\tau^2\phi^{N-1}&(L')\\
\Delta_{g,\phi}W:=\phi^{-N}\div_g(\phi^N\mathring {\mathcal L}_gW)=\frac{1-n}{n}d\tau. &(V')
\end{array}
\right.
$$
 We  	now
make the conformal changes:
$$
\tilde g=\varphi^{N-2}g\;,\;\;\ \tilde \phi=\varphi^{-1}\phi\;,\;\;\ \tilde \sigma=\varphi^{-2}\sigma
\;,\;\;\ \tilde W=\varphi^{N-2}W.
$$
Note that by $(\mathcal C)$, the tensor $\tilde \sigma$ is still TT for $\tilde g$.
Now, we further have :
$$
\tilde \omega:=\omega(\tilde \sigma,\tilde \phi,\tilde W,\tilde g)=\varphi^{-N}\omega
$$
so the scalar equation $(L')$ behaves like $(L)$ as explained  above.  As for the corresponding vectorial equation $(V')$, we now have :
$$
\Delta_{\tilde g,\tilde \phi}{\tilde W}=\Delta_{g,\phi}W,
$$
due to $(\mathcal C)$. In other words, this equation is now  conformally covariant.\\

\noindent{\sc Remarks :}\\
$\bullet$ When $\varphi=\phi$  we find $\Delta_{g,\phi}W=\Delta_{\hat g, 1}\hat W$.\\

\noindent$\bullet$ Since $\omega$ depends on $\phi$, let us give the more explicit form of  $(L')$, namely :
\begin{eqnarray}
\nonumber
\frac{4(n-1)}{n-2}\nabla^*\nabla\phi+R_g\phi-|\sigma|^2\phi^{-N-1}
-2\langle\sigma,\mathring {\mathcal L}_gW\rangle\phi^{-1}-|\mathring {\mathcal L}_gW|^2_g\phi^{N-1}\\
\nonumber
=\frac{1-n}{n}\tau^2\phi^{N-1}. &(L')
\end{eqnarray}

\noindent$\bullet$ Let us consider the operator given by
$$
\mathcal P_g\left(\begin{array}{l}\phi\\ W\\ \sigma\end{array}\right)
:=\left(\begin{array}{l}\phi^{1-N}P_{g,\omega}\phi\\ \Delta_{g,\phi}W\\\phi^{-N}\div_g\sigma\end{array}\right).
$$
This operator is conformally covariant :
$$
\mathcal P_{\tilde g}\left(\begin{array}{l}\tilde\phi\\ \tilde W\\ \tilde\sigma\end{array}\right)
=\mathcal P_g\left(\begin{array}{l}\phi\\ W\\ \sigma\end{array}\right).
$$
The system  $(S')$ (modulo  the third equation added) simply reads :
$$
\mathcal P_g\left(\begin{array}{l}\phi\\ W \\ \sigma \end{array}\right)
=\frac{1-n}n\left(\begin{array}{l}\tau^2\\ d\tau\\ 0\end{array}\right).
$$

\noindent$\bullet$ There has been  other attempts to modify the Lichnerowicz-York method, like
the  ``Conformal Thin-Sandwich''(see for instance   \cite{MaxwellConfetCTS} for a comparison, see also \cite{Maxwellparamconf}). In particular, there  exists already parameterizations where
$$\hat K=\frac{\tau}n\hat g+\phi^{-2}(\sigma+f\mathring {\mathcal L}_g W)$$ but so far the function $f$
would  usually  be given and would not depend on $\phi$, whereas  $f=\phi^N$ is the simplest one to make the system conformally covariant. We will see in section    \ref{autreparam} further possibilities, including conformally covariant ones.\\

\noindent$\bullet$ More generally, for each  conformally covariant differential linear operator of order $k$, we would  interpose a function between the operator and its formal adjoint, in order to produce an operator of order $2k$,  conformally covariant positive  and self adjoint for a weighted measure.

\section{The vectorial equation}\label{equaV}
Let us consider the vectorial laplacian, of Witten type, obtained previously :
$$
\Delta_{g,\phi} W=\phi^{-N}\div_g\phi^N\mathring{\mathcal L}_gW=\div_g\mathring{\mathcal L}_gW-N\mathring{\mathcal L}_gW(\nabla \ln\phi,.).
$$
We choose to work here on a compact manifold.
Since $2\div_g$ is the formal $L^2(d\mu_g)$ adjoint of  $\mathring{\mathcal L}_g$,  we have 
$$
\int_M \langle\Delta_{g,\phi}W,V\rangle_g \phi^{N}d\mu_g=
\frac12\int_M\langle \mathring{\mathcal L}_gW,\mathring{\mathcal L}_gV\rangle_g \phi^{N}d\mu_g
$$
Our Laplacian  is then $L^2(\phi^Nd\mu_g)$ self adjoint and its kernel is reduced  to conformal killing 1-forms.\\

We want to solve for $W$ the equation
$$
\Delta_{g,\phi} W=\frac{1-n}nd\tau\;,\;\;\;(V')
$$
with $g,\phi,\tau$ smooth for simplicity (we could also choose $\phi\in L^{\infty}$).
By the Fredholm alternative, a necessary and sufficient condition for solving $(V')$ is 
orthogonality of the right-hand side $d\tau$ to the kernel:
$$
\int_M \langle X,d\tau\rangle_g\phi^{N}d\mu_g=0,\;\;\;(\perp)
$$
for all $X$ in ker$\;\mathring{\mathcal L}_g$. Uniqueness of $W$ occurs up to the addition of an element of ker$\;\mathring{\mathcal L}_g$.
Note that if  $\tilde X=\varphi^{N-2}X$ ,
$$
\langle \tilde X,Y\rangle_{\tilde g}\tilde\phi^{N}d\mu_{\tilde g}
=\langle X,Y\rangle_g\phi^{N}d\mu_g.
$$
 Hence the condition $(\perp)$ is conformally covariant due to $(\mathcal C)$. It can be written
$$
\int_M \tau d^*_g(\phi^{N}X)d\mu_g=0.\;\;\;(\perp^*)
$$
Of course, if the metric $g$ does not possess some conformal Killing field, the equation
$(V')$ always has a unique solution whatever  the (positive) function
$\phi$ is.\\

If  $\mathcal  K=\ker\mathring{\mathcal L}$ is not trivial,
we may define $\Pi_{g,\phi}$ the orthogonal  projection on
$\mathcal K$ for the weighted scalar product :
$$
(X,Y)_{g,\phi}=\int_M\langle X,Y\rangle_g\phi^{N}d\mu_g.
$$ 
We can then solve
$$
\Delta_{g,\phi} W=\frac{1-n}n(Id-\Pi_{g,\phi})d\tau,\;\;\;(V'_{\Pi})
$$
$W$ being unique if it is chosen  orthogonal to $\mathcal K$ (for the weighted scalar product). \\

\noindent{\sc Remark :} By analogy with the Witten laplacian, we could consider the similar laplacian 
$\phi^{-\frac{N}{2}}\div(\phi^{\frac{N}{2}}\phi^{\frac{N}{2}}\mathring{\mathcal L}\phi^{-\frac{N}{2}})$.
 
\section{Further versions of the method B}\label{methodBbis}
We develop some variants of the method B.
\subsection{A linear version}\label{Line}
Let us linearise  the operator ${\mathcal P}_g$ introduced  in section \ref{paramCC}.
We denote by $(\psi,V,h)$ the variations of $(\phi,W,\sigma)$, 
 $\upsilon$ the variation of $\tau$, we keep here the metric $g$ fixed.
Because of the conformally covariant properties of ${\mathcal P}_g$, the computation of the linearisation can be done at $\phi=1$.
We  find
$$
D{\mathcal P}_g^{(1,W,\sigma)}
\left(\begin{array}{l}\psi\\ V\\ h\end{array}\right)=
\left(\begin{array}{ccc}
P_{g,\sigma,W}& -2\langle \mathring{K}, \mathring{\mathcal L}.\rangle,& 
-2\langle \mathring{K}, .\rangle\\
-N\mathring{\mathcal L}_gW(\nabla_g.,.)&\div_g\mathring{\mathcal L}_g&  0\\
0& 0& \div_g \\
\end{array}\right)
\left(\begin{array}{l}\psi\\ V\\ h\end{array}\right)
$$
where
$$
\mathring{K}=\sigma+\mathring{\mathcal L}_gW,
$$
and
$$
P_{g,\sigma,W}=\frac{4}{n-2}\left((n-1)\nabla^*\nabla-R_g+n|\sigma|^2
+n\langle \sigma, \mathring{\mathcal L}_gW\rangle\right)
$$
$$
\mbox{ }\;\;=\frac{4}{n-2}\left((n-1)\nabla^*\nabla-R_g+\frac n2\left(|\sigma|^2
+|\mathring{K}|^2- |\mathring{\mathcal L}_gW|^2\right)\right).
$$
The latter operator, evaluated at a solution of $(S')$, can be written using the Yamabe laplacian, because
$$
R_g=|\mathring{K}|^2+\frac{1-n}n\tau^2.
$$
The linear version of the system becomes
$$
D{\mathcal P}_g^{(1,W,\sigma)}
\left(\begin{array}{l}\psi\\ V\\ h\end{array}\right)=\frac{1-n}n
\left(\begin{array}{l} 2\tau \upsilon\\ d\upsilon\\ 0\end{array}\right).
$$
If  $\sigma$ is fixed, it forces $h=0$, if not, the global form of the  system suggests  to look at other possibilities.
\subsection{Another conformally covariant parameterization}\label{autreparamCC}
 In the parameterization of the section \ref{paramCC}, we  fixed a  $TT$-tensor $\sigma$ and defined 
$\hat\sigma=\phi^{-2}\sigma$ in order to verify immediately the last equation of the system $(C')$.
Without this particular choice ,  $(C')$ has  $2n+1$ equations. 
Alternatively, we can fix instead a trace free tensor  $\mathring{\sigma}$ and look for the unique
 TT-tensor $\hat{\sigma}$ for $\hat g$,    such that
$$
\phi^{-2}\mathring\sigma=\hat\sigma-\mathring{\mathcal L}_{\hat g} \hat{Y}.
$$
Doing so we should define
$$
\hat{\sigma}=\phi^{-2}(\mathring{\sigma}+\phi^N\mathring{\mathcal L}_g {Y}),\;\;\;\varpi=|\mathring{\sigma}+\phi^N\mathring{\mathcal L}_g {Y}+\phi^N\mathring{\mathcal L}_g {W}|.
$$
It yields  the conformally covariant system
$$
(S'')\left\{\begin{array}{ll} \vspace{-2mm}
\phi^{1-N}P_{g,\varpi}\phi=\frac{1-n}{n}\tau^2&(L'')\\
&\\ \vspace{-2mm}
\Delta_{g,\phi}W=\frac{1-n}{n}d\tau &(V_1'')\\
&\\
\Delta_{g,\phi}Y=-\phi^{-N}\div_g\mathring{\sigma}.&(V_2'')\\
\end{array}
\right.
$$
The matching  linear operator has the form :
$$
\left(\begin{array}{c}
\displaystyle{\phi^{1-N}P_{g,\varpi}\phi}\\
\displaystyle{\Delta_{g,\phi}}\\
\displaystyle{\Delta_{g,\phi}}
\end{array}
\right)'=
\left(\begin{array}{ccc}
P_{g,\mathring \sigma,W,Y}& -2\langle \mathring{K}, \mathring{\mathcal L}.\rangle,& 
-2\langle \mathring{K}, \mathring{\mathcal L}.\rangle\\
-N\mathring{\mathcal L}_gW(\nabla_g.,.)&\Delta_{g,1}&  0\\
-N\mathring{\mathcal L}_gY(\nabla_g.,.)& 0& \Delta_{g,1}\\
\end{array}\right)
$$
where
$$
\mathring{K}=\mathring\sigma+\mathring{\mathcal L}_gY+\mathring{\mathcal L}_gW,
$$
and
$$
P_{g,\mathring \sigma,W,Y}=\frac{4}{n-2}\left((n-1)\nabla^*\nabla-R_g+n|\mathring\sigma|^2
+n\langle \mathring\sigma, \mathring{\mathcal L}_g(W+Y)\rangle\right)
$$
$$
\mbox{ }\;\;=\frac{4}{n-2}\left((n-1)\nabla^*\nabla-R_g+\frac n2\left(|\mathring\sigma|^2
+|\mathring{K}|^2- |\mathring{\mathcal L}_g(W+Y)|^2\right)\right).
$$
This operator evaluated at a solution of
 $(S'')$, can again be written with the Yamabe laplacian, because 
$$
R_g=|\mathring{K}|^2+\frac{1-n}n\tau^2.
$$
Getting back to $(S'')$, we realize that 
 $X=Y+W$ is the only important variable and a natural system to solve is
$$
(S''')\left\{\begin{array}{ll} \vspace{2mm}
\phi^{1-N}P_{g,\varpi}\phi=\frac{1-n}{n}\tau^2&(L''')\\
\Delta_{g,\phi}X=\frac{1-n}{n}d\tau-\phi^{-N}\div_g\mathring\sigma &(V''')\\
\end{array}
\right.
$$
where  the function $\tau$  and the trace free symmetric two tensor $\mathring\sigma$ are given.
The parameterization for the solutions of  $(C)$ will then be
$$
\hat g=\phi^{N-2}g\;,\;\;\hat K=\phi^{-2}(\mathring\sigma+\phi^N\mathring{\mathcal L}_g {X})
\;.\;\;\;(P''')
$$
{\sc Remark} : 
We recover the system $(S')$ if we choose  $\mathring\sigma$ divergence free.\\

\section{Further parameterizations and  matching}\label{autreparam}
Let us propose a general parameterization, from which  the conformal methods A and B, and the conformal thin sandwich method are particular cases.

Let  $\tilde g=\varphi^{N-2}g$ be  another conformal metric.
The York  decomposition $(Y)$, related to $\tilde g$, of the tensor
$$\phi^2\varphi^{-2}( \hat K-\frac{\tau}n\hat g),$$
leads to the parameterization
$$
(P_{\varphi})\;\left\{\begin{array}{l}
\hat g=\phi^{N-2}g\;,\\
 \hat K=\frac{\tau}n\hat g+\phi^{-2}(\sigma+\varphi^N\mathring {\mathcal L}_g W).
\end{array}\right.
$$

Using this parameterization in $(C)$,  and setting
$$
\omega_{\varphi}=\omega(\sigma,\varphi,W,g):=|\sigma+\varphi^N\mathring {\mathcal L}_g W|_g,
$$
we obtain the new system
$$
(S_{\varphi})\left\{\begin{array}{ll}
P_{g,\omega_\varphi}\phi=\frac{1-n}{n}\tau^2\phi^{N-1}&\;\;\;(L_{\varphi})\\
\phi^{-N}\div_g(\varphi^N\mathring {\mathcal L}_gW)=\frac{1-n}{n}d\tau &\;\;\;(V_{\varphi}).
\end{array}
\right.
$$
If  $g$ has no conformal Killing field, we infer from (Y) that,
 for each one form $W$ and each positive function $\varphi$, there exists a unique   1-form $V$ and
a TT-tensor $\sigma'$ such that
$$
\sigma'+\mathring{\mathcal L}_gV=\varphi^N
\mathring{\mathcal L}_{ g} W. \;\;\;\;
$$
It follows that the solutions of the systems  $(S_{\varphi})$
match, but for different  $\sigma$.

Explicitly, $(L_\varphi)$ reads in the following way
\begin{eqnarray}
\nonumber
\frac{4(n-1)}{n-2}\nabla^*\nabla\phi+R_g\phi-|\sigma|^2\phi^{-N-1}
-2\langle\sigma,\mathring {\mathcal L}_gW\rangle\varphi^N\phi^{-N-1}\\
\nonumber-|\mathring {\mathcal L}_gW|^2_g\varphi^{2N}\phi^{-N-1}
=\frac{1-n}{n}\tau^2\phi^{N-1} &(L_\varphi).
\end{eqnarray}

Here, it is important to note that we could let $\varphi$ depend on
$\phi$ and possibly on some other parameters in the (then abusively denoted) system  $(P_{\varphi})$. 

The  conformal method A consists in choosing $\varphi=1$, the conformal method B arises when  $\varphi=\phi$, and for
$\varphi$ a fixed positive function, we obtain  the conformal thin sandwich method.

But many other choices can be made.
For instance $\varphi^N=\phi^NF^2(|\mathring{\mathcal L} W|_g)$, for a given function $F$,  will provide again a conformally covariant system of order of derivation 2.
When possible, the parameterizations $\varphi^N=c|\tau|\phi^N$ or 
$\varphi^N=c\frac{|\tau|}{|\mathring{\mathcal L} W|_g}\phi^N$  seem also to yield  interesting systems.

\section{Comments and prospects}\label{conj}
We give here some comments and hints at a future study of conformally covariant systems, starting with the case of  compact manifolds.

\noindent $\bullet$ {\small\sc Other York decompositions} :
Conformally  covariant or not, in section \ref{autreparam}, some other choices 
of $\varphi=\varphi(\phi,W,\sigma,\tau,t,f,...)$ in $(S_{\varphi})$ may be judicious. \\

\noindent $\bullet$ {\small\sc Constraints with right-hand side} :
The parameterizations works for some other stress energy tensors like  a scalar field
for instance.\\

\noindent $\bullet$ {\small\sc Matching of parameterizations} : 
Due to the matching between the differents parameterizations, already 
existing  solutions of $(S)$ provide  solutions of $(S')$, using a  change of parameters.
If we view the  solutions set of 
 $(C)$ as a manifold, we expect the conformal method to give a local chart,
whereas a conformally covariant method would give a larger chart and allow for a larger choice of
 $(\tau,\sigma)$ (or $(\tau,\mathring\sigma)$ for $(S''')$).\\

\noindent $\bullet$ {\small\sc  Variational study} :
The systems we obtained are not triangular like for the classical conformal method.
However, a variational study seems well suited.\\
 
\noindent $\bullet$ {\small\sc  fixed point} :
As in many previous studies of $(S)$, a classical method  can be used to solve  a system like $(S')$, namely by  considering
a map $T:E\rightarrow E$, with an appropriate function space $E$, defined as follows.
If there is no conformal Killing fields, given
$\psi $ in $E$, we consider the solution  $W$ of
(V') with $\phi$ replaced by $\psi$. We then define $\omega:=\omega(\sigma,\psi,W,g)$
and  solve $(L')$ (see for instance  \cite{MaxwellFreely}).
The solution defines  the map $T(\psi):=\phi$.
Now, we  look for a fixed point  $\phi>0$ of $T$.
If necessary, here, one  could insert some $\phi$ fitting some  $\psi$ in the choice of $\omega$ (and/or in $(L')$) then solve another scalar equation, similar to $(L')$ on $\phi$, for  example linear with respect to $\phi$.
If there exists a conformal Killing field, we can replace in the above process,  the solution  of $(V')$ by that of $(V'_{\Pi})$. The projection disappears in the limit under
appropriate conditions.\\

\noindent $\bullet$ {\small\sc  The limit equation} : In \cite{DGH} a (family of)  limit equation(s)  is proposed for $(S)$. If the limit equation admits only the zero solution, it has the striking  property to guarantee the existence of solutions of $(S)$.
We will have to check if the limit equation measures the asymmetry of the usual parameterization (method A) or if a similar equation exists for the system $(S''')$ (method B).
Depending of the choice of $\varphi$ in section \ref{autreparam},
an associated limit equation  may also appear. 

\section{Appendix}
\subsection{ Vector laplacian}
$\mbox{ }$
$$
\begin{array}{lll}
\div_g\mathring{\mathcal L}_g&=&\nabla^*_g\nabla_g-\Ric_g+\frac{n-2}2dd_g^*\\
&=&\displaystyle{\Delta_{\mbox{\tiny Hodge}}-2\Ric_g+\frac{n-2}2dd_g^*}\\
&=& \displaystyle{d_g^*d+\frac{n}2dd_g^*-2\Ric_g}.
\end{array}
$$

\subsection{ Conformal covariance}\label{cc}
Let us consider three products of  tensor bundles over $M$,
$$
E=E_1\times...\times E_k,\; \;F=F_1\times...\times F_l,\;\;G=G_1\times...\times G_m,
$$
and a differential operator acting on the sections :
$$
P_g:\Gamma(E)\longrightarrow \Gamma(F),
$$
with coefficients  determined 
by $g=(g_1,...,g_m)\in G$. We will say that $P_g$  is conformally covariant
if there exist $a=(a_1,..,a_k)\in\R^k$, $b=(b_1,..,b_l)\in\R^l$ and
$c=(c_1,..,c_m)\in\R^m$ such that for each smooth  section $e$ of $E$, and every smooth function  $\varphi$ on
 $M$, we have
$$
\varphi^b \odot P_{\varphi^c \odot g}(\varphi^a \odot e)=P_g(e),
$$
where  
$$
\varphi^a\odot e=(\varphi^{a_1} e_1,...,\varphi^{a_k} e_k).
$$
A differential system will be said conformally covariant if it can be written  in the form
 $P_g(e)=f$,
for a conformally covariant operator $P_g$.

\def\polhk#1{\setbox0=\hbox{#1}{\ooalign{\hidewidth
  \lower1.5ex\hbox{`}\hidewidth\crcr\unhbox0}}}
  \def\polhk#1{\setbox0=\hbox{#1}{\ooalign{\hidewidth
  \lower1.5ex\hbox{`}\hidewidth\crcr\unhbox0}}} \def\cprime{$'$}
  \def\cprime{$'$} \def\cprime{$'$} \def\cprime{$'$}
\providecommand{\bysame}{\leavevmode\hbox to3em{\hrulefill}\thinspace}
\providecommand{\MR}{\relax\ifhmode\unskip\space\fi MR }
\providecommand{\MRhref}[2]{%
  \href{http://www.ams.org/mathscinet-getitem?mr=#1}{#2}
}
\providecommand{\href}[2]{#2}


\begin{thebibliography}{1}

\bibitem{BartnikIsenberg2004}
Robert Bartnik and Jim Isenberg, \emph{The constraint equations}, The
  {E}instein equations and the large scale behavior of gravitational fields,
  Birkh\"auser, Basel, 2004, pp.~1--38. \MR{2098912}

\bibitem{CookRelatNum2002}
G.~B.~Cook H.~P.~Pfeiffer and S.~A. Teukolsky, \emph{Comparing initial-data
  sets for binary black holes}, Phys. Rev. D (2002), no.~66, 1--17.

\bibitem{IsenbergCMC}
J.~Isenberg, \emph{Constant mean curvature solutions of the {E}instein
  constraint equations on closed manifolds}, Classical and Quantum Gravity
  \textbf{12} (1995), 2249--2273.

\bibitem{Lichnerowicz:44}
A.~Lichnerowicz, \emph{L'int\'egration des \'equations de la gravitation
  relativiste et le probl\`eme des n corps}, Journal de Math\'ematiques Pures
  et Appliqu\'ees \textbf{23} (1944), 37–63.

\bibitem{DGH}
R.~Gicquaud M.~Dahl and E.~Humbert, \emph{A limit equation associated to the
  solvability of the vacuum einstein constraint equations by using the
  conformal method}, Duke Mathematical Journal \textbf{161} (2012), no.~14,
  2669–2697.

\bibitem{MaxwellFreely}
D.~Maxwell, \emph{A class of solutions of the vacuum {E}instein constraint
  equations with freely specified mean curvature}, Math. Res. Lett. \textbf{16}
  (2009), no.~4, 627–645.

\bibitem{Maxwellparamconf}
\bysame, \emph{Initial data in general relativity described by expansion,
  conformal deformation and drift}, arXiv:1407.1467 [gr-qc] (2014).

\bibitem{MaxwellConfetCTS}
\bysame, \emph{The conformal method and the conformal thin-sandwich method are
  the same}, Classical and Quantum Gravity \textbf{31} (2015), no.~14.

\bibitem{York:decompo}
J.~W. York, \emph{Conformally invariant orthogonal decomposition of symmetric
  tensors on riemannian manifolds and the initial value problem of general
  relativity}, Journal of Mathematical Physics \textbf{14} (1973), no.~4,
  456–464.

\end{thebibliography}
\end{document}